\documentclass[a4paper,12pt]{article}

\usepackage{amssymb}
\usepackage{graphicx}
\usepackage{dcolumn}
\usepackage{amsmath,amssymb,mathrsfs}

\begin{document}

\title{Relations between $e$ and $\pi$: Nilakantha's series and Stirling's formula}
\author{V. Yu. Irkhin}
\maketitle
\begin{abstract}
Approximate relations between $e$ and $\pi$ are reviewed, some new connections being established.
Nilakantha's series expansion for $\pi$ is transformed to accelerate its convergence. Its comparison with the standard inverse-factorial expansion for $e$  is performed to demonstrate similarity in several first terms. This comparison, as well as using Dalzell's-like identities, clarifies the origin of the approximate coincidence $e+2\pi \approx 9$. Using Stirling's series enables us to illustrate the relations $\pi^4+\pi^5 \approx e^6$ and $\pi^{9}/e^{8} \approx 10$. The role of Archimedes approximation $\pi=22/7$ and magical 7 number is discussed. A new near-integer relation  $19\pi + 6e = 75.99995…$ is derived using a corresponding lattice construction.


\end{abstract}



The fundamental mathematical constants $e$ and $\pi$ are both transcendental, and their simple combinations are probably transcendental too, except for famous Euler's identity $e^{i\pi}=-1$.

However, there exist very close numerical coincidences in some combinations of $e$ and $\pi$, which seem to be accidental. 
Several examples can be found in \cite{552,Pi}: 
\begin{equation}\pi^2=0.9996(4e-1)\approx 4e-1\end{equation}
\begin{equation}{163}(\pi -e) =68.99966...\end{equation}
\begin{equation}(\pi^4+\pi^5)/ e^6=0.999999956...\,,
\end{equation}
\begin{equation}{ \pi^{9}/e^{8}=9.9998... \approx 10},\end{equation}
\begin{equation} e^\pi-\pi=19.999...,\end{equation}
\begin{equation}{\pi^2 {(\pi -e)^{3/2}}/e=0.9999869...,}\end{equation}
\begin{equation}\pi-e \approx 1- \gamma \end{equation}
with $\gamma \approx 0.57721$ the Euler–Mascheroni constant.
 
For some of such approximate relations,  complicated mathematical justifications were found. 
For example, ``Ramanujan's constant''
\begin{equation}{e^{\pi {\sqrt {163}}}=262\,537\,412\,640\,768\,743.999\,999\,999\,999\,25\ldots \approx 640\,320^{3}+744} \nonumber
\end{equation}
is very close to integer value owing to that 163 is a Heegner number \cite{554}.

Here we try to provide some illustrations and relations for the above near-integer equations.


\section{Nilakantha's series}
First we focus on the simplest relation  
\begin{equation}
	e+2\pi=9.001...,
\end{equation}	
	 or $2(\pi-3) \approx 3 - e$.
We demonstrate a reason of this coincidence by using the series expansions.
For $e$ we use the standard inverse-factorial series which may be represented as
\begin{equation}
	e= \sum_{n=0}^\infty \frac{1}{n!}=3-\frac{1}{3}+\frac{1}{24}+\frac{1}{120}+\frac{1}{720}+ ...
\end{equation}
and is rapidly convergent.

As for  $\pi$, the situation is more complicated. The simple Gregory-Leibniz expansion
\begin{equation}
	\frac{\pi}{4}	= 1 - \frac{1}{3}+\frac{1}{5}-\frac{1}{7}+ ...,
\end{equation}
converges very slowly.
We apply  Nilakantha's series which can be obtained from (10) after acceleration of the convergence:
\begin{eqnarray}
\pi&=& 3+\sum_{n=1}^\infty (-1)^{n+1}  \left( \frac{1}{n+1} + \frac{1}{n}-\frac{4}{2n+1} \right) \nonumber \\
&=&3+\sum_{n=1}^\infty (-1)^{n+1} \frac{1}{n(2n+1)(n+1)}  \nonumber \\
&=&3+\sum_{n=1}^\infty (-1)^{n+1} \frac{4}{(2n+1)^3-(2n+1)} .
\end{eqnarray}%
The series (11) 	was discovered in India already in the 15th–16th century and converges rather rapidly  \cite{55}.

After grouping of pairs of opposite-sign terms, (11) can be transformed as
\begin{eqnarray}
	2\pi&=&6+\frac{1}{3}-\sum_{n=1}^\infty \frac{3}{n(n+1)(4n+1)(4n+3)} \nonumber \\
&=&	6+\frac{1}{3}-\frac{3}{70}-\frac{1}{198}-\frac{1}{780}-...
\end{eqnarray}
so that the convergence becomes still faster.

One can see that the expansions (9) and (12) demonstrate a close similarity of several first terms, which justifies and clarifies the coincidence (8).

The structure of Nilakantha's series also reproduces roughly  the relation (1). We have from (12)
\begin{eqnarray}
	\frac{1}{4}(\pi^2+1) \approx  \frac{1}{4}\left( 3+ \frac{1}{6}- \frac{1}{48}- \frac{1}{200}\right)^2+\frac{1}{4} \nonumber \\
	 \approx 3-\frac{1}{3} +\frac{1}{12}-\frac{1}{32}-0.01 \approx  3-\frac{1}{3} +\frac{1}{24}
\end{eqnarray}
which corresponds to the expansion of $e$ (9).

The relation (8) enables us to present  combinations of $e$ and $\pi$ in terms of one constant and therefore can be useful to clarify other above relations between $e$ and $\pi$. So,  the relation (1) is transformed as
\begin{equation}
	\pi^2+8\pi=35
\end{equation}
which gives
\begin{equation}
	\pi=  \sqrt{51}-4\approx 3.1414...
\end{equation}
Expanding, we obtain
\begin{equation}
	\pi=  7\sqrt{1+2/49}-4\approx 22/7\approx 3.143,
\end{equation}
i.e., Archimedes number.

On substituting (8)  into (2) we have
\begin{equation}
	\pi=  \frac{512}{163}= \frac{2^9}{163}=3.1411...,
\end{equation}
which corresponds to  Stoschek's approximation using powers of two and the  number 163 that is the largest Heegner number \cite{Pi}. 

To comment Eq.(7) we can use the expansion
\begin{equation}
\gamma=  \sum_{n=1}^\infty \frac{|G_n|}{n}=
\frac{1}{2}+\frac{1}{24}+\frac{1}{72}+\frac{19}{2880}+\frac{3}{800}+...,
\end{equation}
where $G_n$ are Gregory coefficients (Bernoulli numbers of the second kind). This expansion is somewhat reminiscent of (9), (12), but converges rather slowly.

\section{Lattice construction}

An analogue of Nilakantha’s approach is given by the integral representation
\begin{equation} 
\pi-3= \int_{0}^{1}\frac{x(1-x)^2}{ 1+x^{2}}\, dx
\end{equation}
which is reduced to (10) by the series expansion. Similarly we have
\begin{equation} 
3- e=\int_{0}^{1}{x(1-x)}e^x, dx.\label{J2}
\end{equation}

A still smaller integral correction is involved into the identity by Dalzell (1944): 
\begin{equation} 
\frac{22}{7}-\pi= \int_{0}^{1}\frac{x^{4}(1-x)^4}{ 1+x^{2}}\, dx= J_1, 
\end{equation}
Also, integrating by part, we obtain from (\ref{J2})
\begin{equation} 
e-\frac{19}{7}=\frac{1}{14}\int_{0}^{1}{x^{2}(1-x)^2} \, dx =J_2/14, 
\end{equation}
Thus both 22/7 and 19/7 are 
integral approximations with the same magical denominator 7, and we have
\begin{equation} 
e + 2\pi - 9 = J_2/14 - 2J_1 = 0.00399 - 0.00253 = 0.00146,
\end{equation}
so that some additional compensation of errors occurs.

For integers $n$, $m$ we have
\begin{equation} n \pi + me = (22 n+19 m)/7 - n J_1 + m J_2/14. \label{jj}
\end{equation}
Since \(22 \equiv 1\) and \(19 \equiv 5\) (mod 7), the number $(22n+19m)/7$ is an integer exactly when \(n + 5m \equiv n - 2m \equiv 0\) (mod 7). Otherwise \(n \pi + me\) is at least about 1/7 away from an integer, for small $n, m$. For \(n - 2m \equiv 0\) (mod 7), defining a lattice, and the error is
$$E(n,m) = -0.001265 n + 0.003996 m.$$
This we can generate approximate equations:
$$
3\pi-2e=4, \, 3e-\pi= 5, \, 5e-4\pi=1 ...
$$
and also for any combinations $n\pi+me$ with integer $n$ and $m$ satisfying $n-2m=7k$.


Generally speaking, the accuracy is  decreasing with  increasing $n. m$.
However, the errors cancel when \(n/m \approx 3.160\). In particular, the $n = 19, m = 6$ point lies on the lattice (19 - 12 = 7). Then \(E = -4.9 \times 10^{-5}\) and we have a beautiful near-integer relation
\(19\pi + 6e = 75.99995...\). 
This is the best lattice point with $m \leq 13$,  76 = 19(22 + 6)/7. 

\section{Stirling's formula}
Another way to relate $e$ and $\pi$  is using Stirling's series \cite{56} which gives a good accuracy even for $n=1$:
\begin{equation}
	\begin{aligned}n!&\simeq 
{\sqrt {2\pi n}}\left({\frac {n}{e}}\right)^{n}\left(1+{\frac {1}{12n}}+{\frac {1}{288n^{2}}}+\cdots
\right)
\nonumber
\end{aligned}
\end{equation}
so that for any $n$
\begin{equation}
e^n\simeq 
{\sqrt {2\pi n}} \frac {n^n}{n!}\left(1+{\frac {1}{12n}}+{\frac {1}{288n^{2}}}+\cdots
\right) \label{19}
\end{equation}

Eq.(\ref{19}) enables one to generate various approximations for $e$ and $\pi$. For example,  we obtain
\begin{equation}
e=\frac{e^{n+1}}{e^{n}}=\left(1+\frac {1}{n}\right)^{n+1/2}\frac {1+\frac {1}{12(n+1)}+\cdots}{1+\frac {1}{12n}+\cdots} .
\end{equation}
For half-integer $n$ in (\ref{19}), using the value of  Euler's Gamma function $\Gamma(1/2)= \sqrt{\pi}$ we derive
\begin{equation}
e^{n+1/2}=\sqrt{2}\frac {(2n+1)^{n+1}}{(2n+1)!!} \left(1+{\frac {1}{6(2n+1)}}+\cdots
\right)
\end{equation}
which yields  a good approximation even for $n=0$: $e=49/18\approx 2.72$.

The same Stirling mechanism gives a general approximation 
\begin{equation}
\pi \simeq e^{2n-1/(6n)}(n!)^2/(2 n^{2n+1}). \label{pipi}
\end{equation}
It gives 3.127, 3.1395 and 3.1410 for $n = 1$, 2, 3.

To demonstrate the usefulness of the relation (\ref{19}), we substitute $e$ from (\ref{19}) with $n=1$ into (3). Expanding up to $1/n^2$ we obtain
\begin{equation}
	\pi^2+\pi= 	
8(1+1/2+1/8+\cdots)\approx 13. \label{22}
\end{equation}
Then we derive again Archimedes number:
\begin{equation}
	\pi=  -1/2+(7/2)\sqrt{1+4/49}\approx 3+1/7 \label{23}
\end{equation}
Using Eq. (\ref{jj}) for $n=1, m=4$ we get	$4e+\pi \simeq 14$. Substituting this into (\ref{pipi}) we obtain Eq.(1).

Another way to derive (\ref{22}) from (3) is to take into account $e^3= 20.08 \approx 20$ and the well-known relation $\pi^3\approx 31$ which follows from the series
\begin{equation}
	\pi^6= 960\sum_{n=0}^\infty \frac{1}{(2n+1)^6}\approx 960. \label{26}
\end{equation}
Then we have
$$
\pi^2+\pi \approx 400/31\approx 13.
$$
Although the accuracy of (15), (\ref{23}) is considerably  worse than that of original relations (1), (3), they demonstrate close connection with ancient investigations of the problem of calculating $\pi$.

Now we provide rough illustrations for other near-integer relations.
To clarify the relation (4) we use (\ref{19}) for both $n=1$ and $n=2$ and derive
\begin{eqnarray}
	e^8&\simeq& 
	64 \pi^3  \left(1+{\frac {1}{12}}+\cdots
	\right)^4 \left(1+{\frac {1}{24}}+\cdots
	\right)^2 \nonumber \\ 
	&\approx& 	\frac {3}{2} 64 \pi^3 = 96\pi^3. \label{27}
\end{eqnarray}
Using (\ref{26}) we see the coincidence of (\ref{27}) and (4). However, the accuracy of this treatment  is  poor in comparison with the real accuracy.


The coincidence (5) can be related to the continued-fraction representation of Gelfond's constant, which also starts from Archimedes approximation~\cite{57}:
\begin{equation}
	e^\pi= [23;7,9,3,1,1,591...]\approx 20+3 + 1/7 +...
\end{equation}
In more detail, we can use $e^3 \approx 20$ to write down
\begin{equation}
	e^\pi= e^3e^{\pi-3}\approx 20(1+1/7+0.01)\approx 20(1+\pi/20)=20+\pi.
\end{equation}
Recently (in 2023), the origin of Eq.(5) has been connected to the sum related to Jacobi theta functions,
\begin{equation}
\sum_{k=1}^\infty(8\pi k^2-2)e^{-\pi k^2}=1.
\end{equation}
Then Eq.(5) is obtained taking into account that $k=1$ term dominates and using Archimedes approximation which results in a compensation of  errors~[1].


The relation (6) is most difficult.  To  obtain a very rough estimation, we take into account (8) and rewrite (6) as 
\begin{equation}
	27 \pi^8 (\pi -3)^{3}/(\pi^2 e)^2\approx 1. \label{31}
\end{equation}
Further on, we use Archimedes $\pi-3 \approx 1/7$, $7^3 \approx 350$, $\pi^2 e \approx 27$ (which follows from the geometric mean corresponding to (8)),  and $$ \pi^8 = 9450 \, \zeta(8) \approx 9450 = 350 \cdot 27$$
with $\zeta(x)$ the Riemann zeta function. Then Eq.(\ref{31}) is approximately  satisfied, but our accuracy is again poor, about 1\% only.

To conclude, we have discussed some near-integer relations of $e$ and $\pi$ numbers.
In some cases, our justification attempts  are only illustrative and not sufficient: surprisingly, the real accuracy of such relations is considerably higher than expected. Therefore, further ideas are required.

Although the subject of the present work can seem somewhat artificial, we note that combinations of
$e$ and $\pi$ occur in physical problems, e.g., in the exact solution of the Kondo model \cite{6}.


\end{document}